\documentclass[bezier,a4paper,leqno,draft]{amsart}
\usepackage[latin1]{inputenc}
\usepackage{amsmath}
\usepackage{latexsym}
\usepackage{amsthm}
\newcommand{\bmath}[1]{\mbox{\mathversion{bold}$#1$}}
\newcommand{\C}{\bmath{C}}
\newcommand{\Z}{\bmath{Z}}
\newcommand{\R}{\bmath{R}}

   \newtheorem{theorem}{Theorem}

   \newtheorem{lemma}[theorem]{Lemma}
 \theoremstyle{definition}
   \newtheorem{definition}{Definition}

 \theoremstyle{remark}

\begin{document}

\centerline{\large \bf Lower bounds for Morse index of constant mean 
curvature tori}

\medskip
\centerline{Wayne Rossman}

\begin{quote} {\bf Abstract:} 
We give three lower bounds for the Morse index of a constant mean curvature 
torus in Euclidean $3$-space in terms of its spectral genus $g$.  The first 
two lower bounds 
grow linearly in $g$ and are stronger for smaller values of $g$, while the 
third grows quadratically in $g$ but is weaker for smaller values of $g$.  
\end{quote}

\section{Introduction}

The Morse index of a complete 
constant mean curvature (CMC) $H$ surface in $\R^3$ with 
$H \neq 0$ is finite if and only if the surface is compact 
\cite{s} \cite{lr}, and is $0$ (i.e. the surface is stable) if and only if 
the surface has genus $0$ and hence is a round sphere \cite{bc}.  It is 
also understood how to find all CMC tori \cite{b2} \cite{ps}.  
Thus, to search for the least possible index of 
unstable CMC surfaces, it is natural to begin with compact CMC tori.  
The simplest of them, the original Wente tori 
\cite{we} \cite{a} \cite{wa} with spectral genus $g=2$, have index 
$\geq 8$ \cite{lnr} \cite{r1} \cite{r2}, suggesting that perhaps no 
unstable CMC surface has index less than $8$.  In this direction, we show 
here that CMC tori with large $g$ must have large index.  (CMC tori 
exist for every $g \geq 2$ \cite{j} \cite{ekt}.)  

\section{Description of CMC $1$ tori}

Any CMC torus can be described as a conformal isometric immersion 
\[ F: \C / \Lambda \to \R^3 \; , \] where $\Lambda$ is a lattice in the 
complex plane $\C$, and the induced Riemannian metric on 
$\C / \Lambda$ is 
\[ ds^2 = e^u \cdot ds^2_{Eucl} \; , \;\;\;\;\;\;\;\; 
\text{where} \;\;\; ds^2_{Eucl}=dx^2 + dy^2 \;  \] is the standard 
Euclidean metric, and $u(z:=x+iy): \C / 
\widetilde{\Lambda} \to \R$ is doubly periodic 
with respect to another lattice $\widetilde{\Lambda}$ of $\C$. 
As CMC tori have no umbilic points \cite{b2}, we may further assume 
the mean curvature and Hopf differential are 
\[ H=1 \;\;\; \mbox{   and   } \;\;\; Q := \langle F_{zz}, \vec{N} \rangle = 
1/2 \; , 
\] where $\vec{N}$ is a unit normal vector to the surface, and hence 
$u$ satisfies the sinh-Gordon equation \[ \partial_z \partial_{\bar{z}} 
u+\sinh u = 0 \; . \] Furthermore, $u$ is smooth, i.e. $u \in C^\infty 
(\C / \widetilde{\Lambda})$.  
(The above facts are explained in more detail in any of \cite{b2} 
\cite{lnr} \cite{r1} \cite{r2} \cite{wa} \cite{we}.)  

Let $\Pi$ (resp. $\widetilde{\Pi}$) represent a fundamental domain of the 
lattice $\Lambda$ (resp. $\widetilde{\Lambda}$).  Suppose that 
$m$ copies of $\widetilde{\Pi}$ translated by vectors in 
$\widetilde{\Lambda}$ can be placed within $\Pi$ with disjoint interiors.  
In other words, we have at least 
$m$ disjoint congruent open regions on the torus 
$F$, each representing a region of double periodicity for $u$.  For CMC tori 
with many symmetries, $m$ can be large; for example, the original Wente tori 
can have arbitrarily large $m$.  Since at the 
very least we may take $\Lambda=\widetilde{\Lambda}$ and $\Pi=
\widetilde{\Pi}$, we may assume \[ m \geq 1 \; . \]

The function $u$ can be described with theta functions (Theorems 7.2 and 
8.1 of \cite{b2}):  
\[ u(z) = 2 \log \left( \frac{\theta(i \mbox{Re}(U z)+D+i \pi (1,1,...,1))}
{\theta(i \mbox{Re}(U z)+D)} \right) \; , \]
where the theta function $\theta$ (as defined in \cite{b2}) 
is determined by a spectral curve of genus $g \geq 2$ defined in \cite{b2}, 
and $D \in i \R^g$ is arbitrary, and $U$ is defined in Theorem 7.1 of 
\cite{b2}.  As the choice of $D$ does not affect the periodicity of the 
surface, $D$ gives a smooth $g-2$ parameter family of CMC tori.  
Furthermore, there is a heirarchy of solutions $v_j: \C / 
\widetilde{\Lambda} \to \R$ to the linearized sinh-Gordon equation 
\begin{equation}\label{operator} 
L(v_j) = 0 \; , \; \; \; \; \; \; L := - \partial_{z}\partial_{\bar{z}} - 
\cosh u \end{equation} given recursively by the 
following procedure \cite{b1}: with the Pauli matrices 
\[ \sigma_1=\begin{pmatrix} 0 & 1 \\ 1 & 0 \end{pmatrix} \; , \;\;\; 
   \sigma_2=\begin{pmatrix} 0 & -i \\ i & 0 \end{pmatrix} \; , \;\;\; 
   \sigma_3=\begin{pmatrix} 1 & 0 \\ 0 & -1 \end{pmatrix} \]
as given in \cite{b2}, define off-diagonal matrices $R_j$ recursively by 
\[ R_1=-\frac{1}{2}u_z\sigma_2 \; , \; \; \; 
R_2=\frac{1}{2}u_{zz}\sigma_1 \; , \; \; \; 
i\left[ R_{k+1},\sigma_3 \right] = -u_z \sum_{n=1}^{k-1} R_n \sigma_1 
R_{k-n} -2 \partial_z R_k \; , \;\;\; k \geq 2 \; . \] Now define $K_j$ 
recursively by 
\[ K_1=-i\sigma_3 \; , \;\;\; K_2=-u_z \sigma_1 \; , \;\;\; 
K_{j+1} = -i \left[ R_j,\sigma_3 \right] - \sum_{i=2}^{j} K_i 
R_{j+1-i} \; , \;\;\; j \geq 2 \; . \] 
For all positive even 
$j \in 2 \cdot \Z^+$, we find that 
$K_{j} = \rho_j \sigma_1$ for some scalar function $\rho_j$.  For example, 
the first three $\rho_j$ are 
\[ \rho_2=-\partial_z u \; , \]
\[ \rho_4=-\frac{1}{2} (\partial_z u)^3+\partial_{z}^{(3)} u \; , \]
\[ \rho_6=-\frac{3}{8} (\partial_z u)^5+\frac{5}{2} \partial_z u 
(\partial_{z}^{(2)} u)^2+\frac{5}{2} (\partial_z u)^2 \partial_{z}^{(3)} u
-\partial_{z}^{(5)} u \; , \] where $\partial_{z}^{(n)}$ represents the 
$n$'th derivative with respect to $z$.  Let 
\begin{equation}\label{Vj} 
v_j := \mbox{Re}(\rho_{j+1}) \mbox{   for }j\mbox{ odd,} \;\;\;\;\; 
v_j := \mbox{Im}(\rho_{j}) \mbox{   for }j\mbox{ even.} \end{equation}
It is proven (with slightly differing notation) in Proposition 3.1 of 
\cite{ps} that these $v_j$ satisfy equation \eqref{operator}.  

As $L$ is elliptic and the $v_j$ are 
defined on the compact space $\C / \widetilde{\Lambda}$, only finitely 
many $v_j$ can be linearly independent; but if the spectral genus of 
the torus is $g$ and the spectral curve is nonsingular in the sense of 
\cite{b2}, then at least the first $g-1$ functions 
$v_1,v_2,...,v_{g-1}$ are linearly independent.  We shall assume that 
the spectral curve is nonsingular, as there is a nonrigorous 
argument in \cite{b2} to show that the singular case never occurs.  

Note that we can also consider the $v_j$ to be defined on $\C / 
\Lambda$ as well as on $\C / 
\widetilde \Lambda$, since $u$ is well-defined on $\C / \Lambda$ as well as 
on $\C / \widetilde \Lambda$.  

\section{Definition of Morse index}

We now turn to the definition of Morse index.  
Let \[ F(t):\C / \Lambda \to \R^3 \; , \;\;\; t\in (-\epsilon, 
\epsilon ) \; , \;\;\; F(0)=F \] be 
a smooth variation of $F$ through immersions $F(t)$.  Let $\vec{E}(t)$ be the 
variation vector field on $F(t)$.  
We can assume, by reparametrizing $F(t)$ for nonzero $t$, 
that $\vec{E}(0)=v 
\vec{N}$, $v \in C^{\infty }(\C / \Lambda)$. Let $a(t)=$area($F(t)$). The 
first variational formula is 
\[ a^{\prime}(0) := \left. \frac{d}{dt} a(t) \right| _{t=0}=
- \int_{\C / \Lambda} v dA \; , 
\]
where $dA = e^u dxdy$. Let $V(t)=$volume($F(t)$), as defined in \cite{bc}.  
Then $V^{\prime }(0)=\int_{\C / \Lambda} v dA$. The variation is 
\emph{volume-preserving} if $\int_{\C / \Lambda} \langle \vec{E}(t),\vec{N}(t) 
\rangle dA(t)=0$ for all $t \in (-\epsilon,\epsilon)$.  In 
particular, $\int_{\C / \Lambda} v dA=0$ when $t=0$, so $a^{\prime}(0)=0$ and 
$F$ is critical for area amongst all volume-preserving variations.

So to see which volume-preserving variations reduce area, one must 
consider which of them make the following 
second variation formula (for volume-preserving variations) negative:  
\begin{equation}\label{2ndvarform} a^{\prime\prime}(0) := 
\left. \frac{d^{2}}{dt^{2}} a(t) \right| _{t=0}=\int_{\C / \Lambda}\{|\nabla
v|^{2}-(4H^{2}-2K)v^{2}\}dA = 4 \int_{\C / \Lambda} vLv dxdy \; , 
\end{equation}
where $K$ and $\nabla$ are the Gaussian curvature and gradient with respect 
to $ds^2$, and $L$ is as in \eqref{operator}.  (Note that actually 
$H=1$ here.)  

\begin{definition}\label{defnofindex}
The {\em index} Ind($F$) is the maximum possible 
dimension of a subspace $\mathcal{U} \subseteq C^\infty(\C / \Lambda)$ for 
which $\int_{\C / \Lambda} v dA = 0$ and $\int_{\C / \Lambda} v Lv dxdy 
< 0$ for all nonzero $v \in \mathcal{U}$.  
\end{definition}

Let $L^2(\C / \Lambda)$ be the Hilbert space of 
measurable functions with finite $L^2$ norm, where 
the standard $L^2$ inner product $\langle \cdot , \cdot \rangle_{L^2}$ on 
$L^2(\C / \Lambda)$ is defined with respect to the metric $ds^2_{Eucl}$.  
The complete set of eigenvalues for $L$ is discrete and can be listed as 
\[ \lambda_1 < \lambda_2 \leq \lambda_3 \leq ... \nearrow +\infty \]
with associated eigenfunctions $\nu_j$, i.e. 
\[ L \nu_j=\lambda_j \nu_j \; , \] where the $\nu_j \in 
C^\infty(\C / \Lambda)$ are chosen to 
form an orthonormal basis for $L^2(\C / \Lambda)$.  

The next lemma is in \cite{bb} and \cite{r1} and other places as well, 
but we include a brief proof. 

\begin{lemma}\label{lemma0} 
\[ \mathcal{K}-1 \leq \mbox{Ind}(F) \leq \mathcal{K} \; , \]
where $\mathcal{K}$ is the number of strictly negative eigenvalues of $L$. 
\end{lemma}

\begin{proof} We have $\lambda_{\mathcal{K}}<0 \leq \lambda_{\mathcal{K}+1}$.  
Since there exist solutions $v_j$ as in \eqref{Vj} that solve 
$L v_j = 0$, in fact $\lambda_{\mathcal{K}+1}=0$.  

Let $\mathcal{U}:=\mbox{span}\{\nu_1,\dots,\nu_\mathcal{K} \}$.  
For any nonzero $\nu=\sum^k_{i=1}a_i\nu_i\in \mathcal{U}$ for 
$a_1,...,a_\mathcal{K} \in \R$, we have
\[ \int_{{\C}/\Gamma}\nu L \nu dxdy=
\sum^\mathcal{K}_{i=1}a^2_i\lambda_i<0 \; . \]  
Therefore, if we choose $\hat{\mathcal U}$ to be a subspace of 
$C^\infty({\C}/\Gamma)$ of maximum 
dimension such that $\int_{{\C}/\Gamma}\hat{\nu} L \hat{\nu} dxdy$ restricted 
to $\hat{\nu} \in \hat{\mathcal U}$ is negative definite, 
then $\dim(\hat{\mathcal U})\geq\dim {\mathcal U}=\mathcal{K}$.  

Suppose that $\dim(\hat{\mathcal U})>\mathcal{K}$, and let 
$P:\mathcal{U}\longrightarrow \hat{\mathcal U}$
be the projection of $\mathcal{U}$ to $\hat{\mathcal U}$ with respect to the 
$L^2$ norm. Since $\dim(P(\mathcal{U}))\leq \mathcal{K}$ there 
exists a $\hat{\nu} \in \hat{\mathcal U}$ with $\langle \hat{\nu},
\hat{\nu} \rangle_{L^2} = 1$ such 
that $\hat{\nu} \perp_{L^2} P(\mathcal{U})$.  It follows that $\hat{\nu} 
\perp_{L^2} \mathcal{U}$ and so 
$\int_{{\C}/\Gamma}\hat{\nu} L \hat{\nu} dxdy \geq0$, a contradiction.  Hence 
$\dim(\hat{\mathcal U})=\mathcal{K}$ and $\mbox{Ind}(F)\leq \mathcal{K}$.  

Now, let $\tau:\hat{\mathcal U} \longrightarrow \R$ be the linear functional
defined by 
\[ \tau(\hat{\nu})=\int_{{\C}/\Gamma} \hat{\nu} dA \; . \]  
Since the dimension of
the kernel of $\tau$ is at least $\mathcal{K}-1$, we have by Definition 
\ref{defnofindex} that 
\[ \mbox{Ind}(F)\geq\dim(\mbox{Ker}(\tau))\geq \mathcal{K}-1 \; . \] 
\end{proof}

\section{Two preliminary lemmas}

Before stating and proving our main theorem, we give two lemmas 
needed for the proof.  The first is a generalization of the 
Euler formula for graphs, which is classical and very well-known, but as we 
will need to allow somewhat nonstandard ``graphs" that include closed loops, 
we give a proof.

Henceforth we shall refer to a compact connected Riemann surface without 
boundary as a {\em closed} Riemann surface.  

\begin{definition}\label{firstdefn}
Let $M$ be a closed Riemann surface.  

1) A finite embedded 
{\em graph-with-loops} \[ \mathcal{G} = \mathcal{G}^\prime \cup 
\sum_{j=1}^r \gamma_j \] is the union of an unoriented 
finite embedded graph $\mathcal{G}^\prime$ in $M$ with 
a finite number of disjoint closed loops $\gamma_1,...,\gamma_r \subset 
M$ that do not intersect $\mathcal{G}^\prime$.  
We allow $\mathcal{G}^\prime$ to be disconnected, and 
we allow $\mathcal{G}^\prime$ to have loop-edges. (Loop-edges are edges whose 
two endpoints are the same vertex, not to be confused with closed loops.)  
There are no vertices on the closed loops $\gamma_j$.  

2) Let $\mathcal{F}$ denote the number of faces of $\mathcal{G}$, that is, 
suppose that $M \setminus 
\mathcal{G}$ consists of $\mathcal{F}$ components, which we call $V_1, ... , 
V_{\mathcal{F}}$.  We must allow the possibility that 
some of the $V_j$ are not homeomorphic to disks, as $\mathcal{G}$ can 
contain closed 
loops.  In fact, some $V_j$ might not even be homeomorphic to planar domains.  
Let the number of 
edges (resp. vertices) of $\mathcal{G}^\prime$ be $\mathcal{E}^\prime$ (resp. 
$\mathcal{V}^\prime$).  Counting each closed loop $\gamma_j$ as one edge, we 
can say that $\mathcal{G}$ has $\mathcal{E}=\mathcal{E}^\prime+r$ edges and 
$\mathcal{V}=\mathcal{V}^\prime$ vertices.  
\end{definition}

\begin{lemma}\label{lemma1}
Let $\mathcal{G} = \mathcal{G}^\prime \cup 
\sum_{j=1}^r \gamma_j$ be a graph-with-loops on a closed Riemann 
surface $M$.  Suppose 
that $\mathcal{G}^\prime$ is not empty (i.e. $\mathcal{G}^\prime$ has at 
least one vertex).  Let 
$\chi (M)=2-2 \cdot \text{genus}(M)$ be the 
Euler characteristic of $M$.  Then, with $\mathcal{F}$ and 
$\mathcal{E}$ and $\mathcal{V}$ as in part 2 of 
Definition \ref{firstdefn}, we have 
\begin{equation}\label{euler} 
\mathcal{F}-\mathcal{E}+\mathcal{V} \geq \chi (M) \; . \end{equation}
\end{lemma}

{\bf Remark.} 
Strict inequality can occur in Equation \eqref{euler}.  
As a simple example, consider a graph on a torus that has one loop-edge 
$e$ and one vertex $p$, where $e$ lies in a homotopically trivial loop and 
both its ends connect to $p$.  Then strict inequality will hold.  This 
example is too simple to occur in the proof of Theorem 
\ref{thm1}, but it illustrates why we can only invoke the inequality 
\eqref{euler} (and cannot assume equality) in that proof.  

{\bf Remark.} 
Lemma \ref{lemma1} does not hold without the assumption that 
$\mathcal{G}^\prime$ is nonempty, and a simple 
counterexample is to let $\mathcal{G}$ consist of only a single closed 
loop $\gamma_1$ in the sphere $S^2$.  

{\bf Remark.} 
Equality can hold in Equation \eqref{euler} even if some components $V_j$ of 
$M \setminus \mathcal{G}$ are not 
homeomorphic to disks.  For example, consider a graph 
$\mathcal{G} = \{ e \} \cup \{ p \} \cup \gamma_1$ on the sphere $S^2$ that 
has one loop-edge $e$ and one vertex $p$, where both ends of $e$ connect to 
$p$, and also includes a single closed loop $\gamma_1$ disjoint from 
$e$ and $p$.  Then equality holds in Equation \eqref{euler}, even though 
one of the components of $S^2 \setminus \mathcal{G}$ is homeomorphic to an 
annulus.  

However, if $\mathcal{G}$ does not contain any closed loops, i.e. if 
$\mathcal{G}=\mathcal{G}^\prime$, 
then equality holds in Equation \eqref{euler} if and only if 
each component of $M \setminus \mathcal{G}$ is homeomorphic to a disk, as 
follows from the 
generalized Euler formula (Equation \eqref{genEuler} below).  

\begin{proof}
As $\mathcal{G}$ may contain closed loops, it is not a graph in the usual 
sense, so we cannot immediately apply the generalized Euler formula.  We 
will add edges and vertices to 
$\mathcal{G}$ until it becomes a graph in the usual sense, and then 
apply the formula.  

Let $V_1,...,V_{\mathcal{F}}$ be the faces of $\mathcal{G}$ as in 
part 2 of Definition 
\ref{firstdefn}.  
Let $\mathcal{R}$ be the union of open regions $V_j$ such that the boundary 
$\partial V_j$ has nonempty intersection with $\mathcal{G}^\prime$.  
Note that $\mathcal{R}$ is not empty, because $\mathcal{G}^\prime$ is not 
empty.  If 
$\mathcal{R}=M \setminus \mathcal{G}^\prime$, then 
$\mathcal{G}$ contains no closed loops (i.e. $r=0$) and 
$\mathcal{G}=\mathcal{G}^\prime$ is a graph in the standard sense.  If 
$\mathcal{R} \neq M \setminus \mathcal{G}^\prime$, then 
there must exist some $j_0$ such that $V_{j_0} \subseteq \mathcal{R}$ and 
$\partial V_{j_0}$ contains a loop $\gamma_{j_1}$ for some 
$j_1$.  By reordering the $\gamma_j$ if necessary, we may assume $j_1=r$.  
Since $V_{j_0}$ has 
boundary components in both $\mathcal{G}^\prime$ and $\gamma_r$, it is 
not simply-connected and we can add a vertex $p$ at 
any place along $\gamma_r$ and 
connect $p$ by an edge $e$ to some vertex of $\mathcal{G}^\prime$ so that the 
interior of $e$ lies in $V_{j_0}$ and $e$ does not disconnect 
$V_{j_0}$.  Including $p$ and $e$ results in a graph-with-loops 
\[ \mathcal{G}_{r-1} = 
\mathcal{G}_{r-1}^\prime \cup \sum_{j=1}^{r-1} \gamma_j \; , 
\;\;\;\;\; \text{where} \;\;\; \mathcal{G}_{r-1}^\prime = \mathcal{G}^\prime 
\cup \{ e \} \cup \{ p \} \cup \{ \gamma_r \setminus \{ p\}  \}\]  
and the closed loop 
$\gamma_r$ has become the loop-edge $\gamma_r \setminus \{ p \}$ in 
$\mathcal{G}_{r-1}^\prime$.  In particular, $\mathcal{G} \subseteq 
\mathcal{G}_{r-1}$ (as sets in $M$).  (The subscript $r-1$ in 
$\mathcal{G}_{r-1}$ indicates that $\mathcal{G}_{r-1}$ has $r-1$ closed 
loops.) 

Denoting 
by $\mathcal{F}_{r-1}$, $\mathcal{E}_{r-1}$ and $\mathcal{V}_{r-1}$ the 
number of faces, edges and vertices of $\mathcal{G}_{r-1}$, we have that 
$\mathcal{F}_{r-1}=\mathcal{F}$, $\mathcal{E}_{r-1}=\mathcal{E}+1$ and 
$\mathcal{V}_{r-1}=\mathcal{V}+1$, hence 
\[ \mathcal{F}_{r-1}-\mathcal{E}_{r-1}+\mathcal{V}_{r-1} = 
\mathcal{F}-\mathcal{E}+\mathcal{V} \; . \]
Repeating this procedure $r-1$ more times, we can make a sequence of 
graphs-with-loops $\mathcal{G}_{r-1}$, $\mathcal{G}_{r-2}$, ..., 
$\mathcal{G}_{0}$ that change all the closed 
loops $\gamma_j$ one by one into loop-edges.  Each $\mathcal{G}_{s}$ has 
$s$ closed loops, and $\mathcal{G}_s \subseteq \mathcal{G}_t$ (as sets in $M$) 
when $t \leq s$.  The final graph-with-loops 
$\mathcal{G}_0$ has no closed loops and hence is actually a graph 
in the standard sense, i.e. $\mathcal{G}_0 = \mathcal{G}_0^\prime$.  We have 
$\mathcal{G} \subseteq \mathcal{G}_0$ (as sets in $M$), and the 
number of faces $\mathcal{F}_0$ of $\mathcal{G}_0$ equals $\mathcal{F}$.  
Furthermore, defining $\mathcal{E}_0$ and $\mathcal{V}_0$ as the 
number of edges and vertices of $\mathcal{G}_0$, we have 
\[ \mathcal{F}_0-\mathcal{E}_0+\mathcal{V}_0 = 
\mathcal{F}-\mathcal{E}+\mathcal{V} \; . \]  So it is sufficient to show that 
\begin{equation}\label{finalgoal}
\mathcal{F}_0-\mathcal{E}_0+\mathcal{V}_0 \geq \chi (M) \; . 
\end{equation}  

The graph $\mathcal{G}_0$ will have loop-edges if $r \geq 1$, but it 
has no closed loops, so the 
generalized Euler formula (see, for example, Chapter 9 of \cite{g}) 
can be applied to $\mathcal{G}_0$.  Letting $V_{1,0},...,
V_{\mathcal{F},0}$ be the faces of 
$\mathcal{G}_0$, define $\chi (V_{j,0})$ to be the Euler 
characteristic of $V_{j,0}$.  ($\chi (V_{j,0})$ can 
be computed using any true triangulation of $V_{j,0}$.)  The 
generalized Euler formula says that 
\begin{equation}\label{genEuler} 
\mathcal{V}_0-\mathcal{E}_0+\sum_{j=1}^{\mathcal{F}} 
\chi (V_{j,0}) = \chi (M) \; . \end{equation}  
Since $\chi (V_{j,0}) \leq 1$, this implies Equation \eqref{finalgoal}.  
\end{proof}

The next lemma is the Courant nodal domain theorem.  
The proof is well known (see \cite{c} or \cite{ch}, for 
example), but we include it here because 
we add a potential function to the Laplacian operator (this has little effect 
on the proof), and also because we consider the case of multiple 
eigenvalues.  

Let $M$ be a closed Riemann surface with smooth metric $ds^2$.  
Let $dA$ and $\nabla$ and $\triangle$ be the area form and gradient and 
Laplace-Beltrami operator on 
$M$ associated to $ds^2$.  We choose the sign of $\triangle$ so that 
$\int_M \phi \triangle \phi dA = + \int_M |\nabla \phi|^2 dA$ for general 
smooth functions $\phi$ on $M$.  Consider the operator 
\begin{equation}\label{generaleqn} 
c \cdot \triangle + V 
\end{equation} 
on $M$, where $c$ is a positive constant and $V$ is a smooth bounded function 
on $M$.  We write the complete set of eigenvalues for $c \cdot 
\triangle + V$ as 
\[ \lambda_1 < \lambda_2 \leq \lambda_3 \leq ... \nearrow +\infty \]
with associated smooth eigenfunctions $\nu_j$, i.e. 
\[ (c \cdot \triangle + V) \nu_j=\lambda_j \nu_j \; , \] where the $\nu_j$ are 
chosen to form an orthonormal basis for the function space $L^2(M)$ on $M$.  

The Sobolev space $H^1(M)$ of $M$ is defined to be 
functions in $L^2(M)$ whose weak first derivatives exist and are also in 
$L^2(M)$.  (Here the $L^2$ norm is defined with respect to the metric 
$ds^2$ on $M$.)  There is a standard $H^1$ norm, with respect to 
$ds^2$, which makes $H^1(M)$ a Hilbert space.  

\begin{definition}
For a function $\nu : M \to \R$, the set $\nu^{-1}(0)$ is 
the {\em nodal set} of $\nu$, and each component of 
$M \setminus \nu^{-1}(0)$ is a {\em nodal domain} of $\nu$.  
\end{definition}

\begin{lemma}\label{courant}
The number of nodal domains of $\nu_j$ is at most $j$, for 
every $j \in \Z^+$.  Furthermore, in the case of a multiple eigenvalue 
$\lambda_i=\lambda_{i+1}=...=\lambda_{i+k}$, the number of nodal domains 
of any eigenfunction $\nu \in \text{span}\{\nu_i,...,\nu_{i+k}\}$ is at most 
$i$.  
\end{lemma}

\begin{proof}
Assume $\nu_j$ has at least $j+1$ nodal domains, $j+1$ of which are 
$\Omega_1,...,\Omega_{j+1}$.  Define the
function $\phi_k$ by 
\[ \phi_k=\nu_j \text{ on }\Omega_k\; , \;\;\;\;\; \phi_k = 0 
\text{ elsewhere, } \] 
for $k=1,...,j$ (we exclude $k=j+1$).  By Theorems 
2.2 and 2.5 in \cite{c}, the boundary of 
each $\Omega_k$ is piece-wise smooth, and consists of a finite number of 
smooth curves of finite length and finite total curvature, so the 
weak first derivatives of $\phi_k$ exist and are bounded.  Hence 
$\phi_k \in H^1(M)$.  
Clearly, $\langle \phi_{k_1},\phi_{k_2} \rangle_{L^2}=0$ for any unequal 
$k_1$ and $k_2$, so $\text{span}\{\phi_1,...,\phi_{j}\}$ is of dimension 
$j$.  
Then, since $\text{span}\{\nu_1,...,\nu_{j-1}\}$ is of dimension $j-1$, 
there exists some linear combination 
\[ \phi=\sum_{k=1}^j a_k \phi_k \; , \;\;\; a_k \in \R \] that is 
$L^2$-perpendicular to span$\{\nu_1,...,\nu_{j-1}\}$, i.e. 
\begin{equation}\label{proofone}
\phi \in \left( \text{span}\{\nu_1,...,\nu_{j-1}\} \right)^{\perp_{L^2}} \; . 
\end{equation} 
Furthermore, 
\begin{equation}\label{prooftwo}
\phi \in H^1(M) \; . 
\end{equation} 
Because $\phi$ is a linear combination of the $\phi_k$, the 
Rayleigh quotient $R(\phi)$ of $\phi$ for the operator 
$c \cdot \triangle + V$ satisfies 
\begin{equation}\label{proofthree}
R(\phi) := 
\frac{\int_M \phi ((c \cdot \triangle+V)\phi) dA}{\int_M \phi^2 dA} = 
\lambda_j \; . 
\end{equation}
By \eqref{proofone}, \eqref{prooftwo} and \eqref{proofthree}, it follows 
that $\phi$ is an eigenfunction with 
eigenvalue $\lambda_j$.  (For arguments that show this, 
see \cite{bu}, \cite{be}, \cite{rs} or \cite{u}, for example.  In \cite{be}, 
the argument is given in full detail for the Laplacian operator on 
compact Riemannian manifolds, and the same argument can be applied to the 
operator $c \cdot \triangle + V$ here.)  

However, the eigenfunction $\phi$ is 
identically zero on $\Omega_{k+1}$.  This contradicts the maximum 
principle (see \cite{pw}, for example), and proves the first sentence of the 
lemma.  

Now suppose that $\lambda_i=\lambda_{i+1}=...=\lambda_{i+k}$ and 
that $\nu$ is any function in $\text{span}\{\nu_i,...,\nu_{i+k}\}$.  We 
are free to choose the functions $\nu_1,\nu_2,\nu_3,...$ so that 
$\nu / || \nu ||_{L^2} = \nu_i$, and then the argument in the previous two 
paragraphs shows that $\nu$ has at most $i$ nodal 
domains.  Since $\nu \in \text{span}\{\nu_i,...,\nu_{i+k}\}$ is arbitrary, the 
second sentence of the lemma is shown.  
\end{proof}

\section{Linear lower bounds for the Morse index of the tori}

\begin{theorem}\label{thm1}
If the torus $F$ has spectral genus $g$ and $m \geq 1$ disjoint congruent open 
pieces representing regions of double periodicity for $u$, then 
\[ \mbox{Ind}(F) \geq m \cdot \left[ \frac{g-1}{3} \right] - 2 \; , \] where 
$\left[ r \right]$ denotes the greatest integer less than or equal to a 
real number $r$.  
\end{theorem}

\begin{proof}
The theorem is vacuously true if $g \leq 3$, so we may assume $g \geq 4$.  
Let the $v_j$ be as in Equations \eqref{operator} and \eqref{Vj}.  

Since $ \{ \sum_{j=1}^{g-1} a_j v_j \, | \, a_j \in \R \} $ is a 
$g-1$ dimensional space in the null-space of 
the operator $L$, we can choose $a_j$ so 
that $v=\sum_{j=1}^{g-1} a_j v_j$ has zeroes of order $1$ at $\left[ 
\frac{g-1}{3} \right]$ arbitrary distinct points 
$p_j \in \C / \widetilde{\Lambda}$ for 
$j=1,...,\left[ \frac{g-1}{3} \right]$; that is, 
at each of these points $p_j$ we have \[ v(p_j)=\partial_x v(p_j)=
\partial_y v(p_j)=0 \; . \]  Thus at 
each of these points $p_j$ the nodal set (zero set) of $v$ is locally 
a crossing with equiangular intersection of at least two curves, 
by Theorem 2.5 of \cite{c}.  

Since we may consider the $v_j$ to be functions 
on $\C / \Lambda$, we now consider $v$ to be a function defined on 
$\C / \Lambda$.  Let $\mathcal{G}$ be the nodal set of $v$ on 
$\C / \Lambda$.  Theorems 2.2 and 2.5 in \cite{c} imply that 
$\mathcal{G}$ forms a graph-with-loops on $\C / \Lambda$ with smooth edges and 
isolated vertices where an even number of edges meet equiangularly.  
We also note that 
\begin{equation}\label{line0pt5}
\# (\mbox{vertices of }\mathcal{G}) \geq m \cdot 
\left[ \frac{g-1}{3} \right] \geq 1 \; . 
\end{equation}
Lemma \ref{lemma1} and the first remark following it imply 
\begin{equation}\label{line1} 
\# (\mbox{components of }(\C / \Lambda) \setminus \mathcal{G}) \geq 
\# (\mbox{edges of }\mathcal{G}) - \# (\mbox{vertices of }\mathcal{G}) 
\; . \end{equation} 
Since each vertex of $\mathcal{G}$ has degree at least four, 
we have the inequality 
\begin{equation}\label{line2} 
\# (\mbox{edges of }\mathcal{G}) \geq 2 
(\# (\mbox{vertices of }\mathcal{G})) \; . \end{equation} Equations 
\eqref{line0pt5}, \eqref{line1} and \eqref{line2} combine to give 
\begin{equation}\label{line3} 
\# (\mbox{components of }(\C / \Lambda) \setminus \mathcal{G}) 
\geq m \cdot \left[ \frac{g-1}{3} \right] \; . \end{equation}

Let $\lambda_j$ and $\nu_j$ be the eigenvalues and eigenfunctions 
(as defined in Section 3) of $L$.  Since $L$ has 
$\mathcal{K}$ negative eigenvalues, and $0$ is a multiple eigenvalue of 
order at least $g-1$, we have 
$\lambda_{\mathcal{K}} < 0$ and 
$0=\lambda_{\mathcal{K}+1}=...=\lambda_{\mathcal{K}+k}$ and 
$0<\lambda_{\mathcal{K}+k+1}$ for some 
$k \geq g-1$.  So we have 
\[ v \in \text{span}\{v_{1},...,v_{g-1}\} 
\subseteq \text{span}\{\nu_{\mathcal{K}+1},...,\nu_{\mathcal{K}+k}\} \; . \]
Consider $\C / \Lambda$ with the Euclidean metric $ds^2_{Eucl}$, then 
the associated Laplace-Beltrami operator is $\triangle_{Eucl}=
-\partial_x \partial_x - \partial_y \partial_y$.  Furthermore, 
$L = (1/4) \cdot \triangle_{Eucl} - \cosh u$ is of the form in 
\eqref{generaleqn}, and then Lemma \ref{courant} implies 
that $v$ has at most $\mathcal{K}+1$ nodal domains, so 
\begin{equation}\label{line4} \mathcal{K} - 1 \geq 
\# (\mbox{components of }(\C / \Lambda) \setminus \mathcal{G}) -2 \; . 
\end{equation} 
Combining Equations \eqref{line3} and \eqref{line4} with 
Lemma \ref{lemma0}, we have 
\[ \mbox{Ind}(F) \geq \mathcal{K} - 1 \geq 
\# (\mbox{components of }(\C / \Lambda) \setminus \mathcal{G}) -2 \geq 
m \cdot \left[ \frac{g-1}{3} \right] - 2 \; . \] 
\end{proof}

This result can be improved if $i \mbox{Re}(U z_0)+D =(0,0,...,0) 
= \vec{0}$ for 
some $z_0$.  In this case we can translate the parameter $z \to z+z_0$ and 
assume that $D=\vec{0}$.  
The theta function satisfies the symmetry $\theta (\omega)=\theta(-\omega)$, 
so when $D=\vec{0}$, it follows that $u(z)=u(-z)$.  Let $0,w_1,w_2,w_3$ be 
the four distinct points in $\C / \widetilde \Lambda$ such that 
$ 2 \cdot w_l$ is contained in the lattice $\widetilde \Lambda$ for 
$l=1,2,3$.  The 
fact that $u(z)=u(-z)$ implies that also $u(w_l+z)=
u(w_l-z)$ for any of the three half-periods $w_l,l=1,2,3$.  
From the recursions defining $v_j$ and leading up to \eqref{Vj}, we have 
\begin{equation}\label{antisymmetry} 
v_j(z)=-v_j(-z)\mbox{  and  }v_j(w_l+z)=-v_j(w_l-z) \; . 
\end{equation}

\begin{theorem}\label{thm2}
With the same conditions as in Theorem 
\ref{thm1}, if we also have $D=\vec{0}$, then 
\[ \mbox{Ind}(F) \geq m \cdot \left( \left[ \frac{g-1}{3} \right] + \left[ 
\min(\frac{g-1}{3},4) \right] \right) - 2 \; . \] 
\end{theorem}

\begin{proof}
As in the proof of Theorem \ref{thm1}, we can choose the $p_j 
\in \C / \widetilde{\Lambda}$ arbitrarily.  Furthermore, we define 
$v=\sum_{j=1}^{g-1} a_j v_j$ just as in that proof, so that $L v=0$ and 
$v$ has zeroes of order $1$ at the 
$p_j$ for $j=1,...,\left[ \frac{g-1}{3} \right]$.  

Choosing the first four points $p_j$ to be $p_1=0$, $p_2=w_1$, $p_3=w_2$, 
$p_4=w_3$, then the antisymmetry 
\eqref{antisymmetry} implies $v(p_j+z)=-v(p_j-z)$ for $j \leq 4$, hence 
the nodal set of $v$ locally has $2 \ell$ edges intersecting at $p_j$ with 
$\ell$ odd - in particular, 
there are at least six edges intersecting at $p_j$, for $j \leq 4$.  The 
result then follows exactly as in the proof of Theorem \ref{thm1}, simply 
by noting that in this case one can add \[ m \cdot \left[ 
\min(\frac{g-1}{3},4) \right] \] to the right-hand side of 
Equation \eqref{line2} and this equation will still hold (because the 
$p_j$ for $j \leq 4$ each have at least six adjacent edges).  
\end{proof}

\section{A quadratic lower bound for the Morse index of the tori}

We conclude with another method for finding lower bounds for the index of 
closed CMC tori.  This second method gives estimates that are weaker by many 
orders of magnitude for smaller values of $g$, but it has the advantage that 
its estimates grow quadratically in $g$.  

Let $M$ be a closed Riemann surface of genus $G$ with smooth metric 
$ds^2$ that is conformal to the complex structure of $M$.  Let $dA$ and 
$\triangle$ be the area form and Laplace-Beltrami operator on $M$ with respect 
to $ds^2$, with the same sign convention for $\triangle$ as in Section 4.  
Then, let 
\[ \beta_1 < \beta_2 \leq \beta_3 \leq ... \] be the complete set of 
eigenvalues of $\triangle$.  
(Each $\beta_k$ is repeated the number of times equal to its 
multiplicity.)  Let \[ 
A=\int_M dA \] be the area of $M$.  Theorem 0.5 in \cite{k} tells us that 
there exists a universal constant 
$\tilde{C} > 0$ so that for all $k \geq 1$, 
\begin{equation}\label{korev} 
\beta_k \leq \tilde{C} (G+1) \frac{k}{A} \; . \end{equation} 
We wish to apply Equation \eqref{korev} to the CMC 1 isometric immersions 
$F: \C / \Lambda \to \R^3$ in Section 2, with any spectral genus 
$g \geq 2$.  So we take 
$M=\C / \Lambda$ and hence $G=1$, and we take $ds^2=e^u \cdot ds^2_{Eucl}$ on 
$\C / \Lambda$. 
Let $K$ and $H=1$ be the Gauss and mean curvatures of $F(\C / \Lambda)$, 
considered as functions on $\C / \Lambda$.  

We have the following variational characterizations for the $k$'th 
eigenvalues $\beta_k-2$ and $\hat{\beta}_k$ of the operators 
$\triangle-2$ and $\triangle-4 H^2+2 K$:  
\begin{equation}\label{vc1}
\beta_k-2 = \inf_{M_{k}} 
\left( \sup_{\psi \in M_{k}, 
\psi \neq 0} \frac{\int_{\C/\Lambda} 
\psi ((\triangle - 2) \psi) dA}{\int_{\C/\Lambda} \psi^2 dA} \right) \; , 
\end{equation}
\begin{equation}\label{vc2}
\hat{\beta}_k = \inf_{M_{k}} 
\left( \sup_{\psi \in M_{k}, 
\psi \neq 0} \frac{\int_{\C/\Lambda} 
\psi ((\triangle-4 H^2+2 K) \psi) dA}{\int_{\C/\Lambda} \psi^2 dA} \right) \; 
, \end{equation} 
where $M_{k}$ runs through all $k$ dimensional subspaces of 
$C^\infty(M)$.  (For arguments that show this, 
see \cite{bu}, \cite{be} or \cite{u}, for example.  Again we remark that 
the argument is given in full detail for the Laplacian operator in \cite{be}, 
and that same argument can be applied to the operators here).  

Noting that $-4 H^2 + 2 K = 2K - 4 \leq -2$, 
the variational characterizations \eqref{vc1} and \eqref{vc2} imply 
\begin{equation}\label{twoEVs}
\beta_k-2 \geq \hat \beta_k \; . 
\end{equation}
For the operator $L$, as in Equation 
\eqref{operator}, and for the eigenvalues $\lambda_k$ of $L$, as in 
Section 3, we have the following variational characterization: 
\[
\lambda_k = \inf_{M_{k}} 
\left( \sup_{\psi \in M_{k}, 
\psi \neq 0} \frac{\int_{\C/\Lambda} 
\psi L \psi dxdy}{\int_{\C/\Lambda} \psi^2 dxdy} \right) \; . 
\]
Then, since the final equality of Equation 
\eqref{2ndvarform} holds for any smooth function $v$, we have 
\begin{equation}\label{vc3}
\lambda_k = \inf_{M_{k}} 
\left( \sup_{\psi \in M_{k}, 
\psi \neq 0} \frac{\int_{\C/\Lambda} 
\psi ((\triangle-4 H^2+2 K) \psi) dA}{4 \int_{\C/\Lambda} \psi^2 dxdy} \right) 
\; , \end{equation} 
Since the respective forms $dA$ and $dxdy$ differ by a positive factor 
$e^u$ bounded away from both $0$ and $\infty$, the variational 
characterizations \eqref{vc2} and \eqref{vc3} imply that the $k$'th 
eigenvalue $\hat \beta_k$ is negative if and only if the $k$'th eigenvalue 
$\lambda_k$ is negative.  Hence, by Lemma \ref{lemma0}, Ind$(F)$ is 
greater than or equal to one less than the number of negative eigenvalues 
$\hat \beta_k$ of the operator $\triangle-4 H^2+2 K$.  Then, with 
\[ A=\text{area}(F(\C / \Lambda)) \; , \] 
the inequalities \eqref{korev} and \eqref{twoEVs} imply 
\begin{equation}\label{korevaarway}
\mbox{Ind}(F) \geq 
\# \{k \; | \; \beta_k < 2 \} - 1 \geq 
\left[ \frac{2 A}{\tilde{C} (G+1)} \right] - 2 = 
\left[ \frac{A}{\tilde{C}} \right] - 2 \; . 
\end{equation}
When the CMC 1 torus $F(\C / \Lambda)$ has spectral genus $g$, 
it is shown in \cite{flpp} that 
\begin{equation}\label{flppway}
A \geq \frac{\pi}{4} \left( (g+2)^2-\frac{1}{2} (1+(-1)^g) \right) \; . 
\end{equation}
Combining inequalities \eqref{korevaarway} and \eqref{flppway} results in: 

\begin{theorem}\label{thm3}
There exists a universal constant $C > 0$ such that if the torus $F$ has 
spectral genus $g$, then 
\[ \mbox{Ind}(F) \geq C \left( (g+2)^2-\frac{1}{2} (1+(-1)^g) \right) -2 
\geq C g^2 - 2 \; . 
\] 
\end{theorem}

Because $\tilde{C}$ in \cite{k} is on the order of $10^7$, $C$ is on the order 
of $10^{-7}$, so even though the lower bound in Theorem \ref{thm3} grows 
quadratically in $g$, it is much weaker than the lower bounds in Theorems 
\ref{thm1} and \ref{thm2} for smaller values of $g$.  

{\bf Remark.} 
One could make similar arguments using Theorem 16 of \cite{ly} instead of 
Theorem 0.5 in \cite{k}, but 
then one has a lower bound that also depends on the diameter and a lower 
bound for the Gaussian curvature of the surface.

\par
\setlength{\parindent}{0in}
\footnotesize
\vspace{0.2cm}
{\sc Wayne Rossman}:\, Department of Mathematics, Faculty of Science, 
Kobe University, 
Rokko, Kobe 657-8501, Japan. \, {\it E-mail}: 
wayne@math.kobe-u.ac.jp \, {\it web page}: 
www.math.kobe-u.ac.jp \par
\par


\begin{thebibliography}{20}
\bibitem{a}
   U. Abresch, 
   {\itshape Constant mean curvature tori in terms of elliptic functions},
   J. reine u. angew Math. 374, 169-192 (1987).
\bibitem{bu} 
     S. Bando and H. Urakawa, 
     {\itshape Generic properties of the eigenvalue of the 
     Laplacian for compact Riemannian manifolds}, Tohoku 
     Math. Journal 35, 155-172 (1983).  
\bibitem{bb} 
     L. Barbosa, P. B\'{e}rard, 
     {\itshape A ``twisted" eigenvalue problem and applications to geometry}, 
     Journal Math. Pures Appl. 79, 427-450 (2000).  
\bibitem{bc}
    L. Barbosa and M. do Carmo, 
    {\itshape Stability of hypersurfaces with constant mean curvature},  
     Math. Z. 185, 339-353 (1984).
\bibitem{be}
    P. B\'{e}rard, 
    {\itshape Spectral Geometry: Direct and Inverse 
    Problems},  Lecture Notes in Math. 1207, Springer-Verlag (1986).  
\bibitem{b2}
   A. I. Bobenko,
   {\itshape Constant mean curvature surfaces and integrable equations}, 
   Russian Math. Surveys 46:4, 1-45 (1991).
\bibitem{b1}
   A. I. Bobenko,
   personal communication.
\bibitem{ch}
   I. Chavel,
   {\itshape Eigenvalues in Riemannian Geometry}, 
   Academic Press Inc., 1976.
\bibitem{c}
   S.-Y. Cheng,
   {\itshape Eigenfunctions and nodal sets}, 
   Comment. Math. Helvetici 51, 43-55 (1976).
\bibitem{ekt}
   N. M. Ercolani, H. Knorrer and E. Trubowitz, 
   {\itshape Hyperelliptic curves that generate constant mean curvature 
   tori in $\R^3$}, 
   Integrable systems: The Verdier Memorial conference, Progress in 
   Math. 115, Birkhauser 1993.
\bibitem{flpp}
   D. Ferus, K. Leschke, F. Pedit and U. Pinkall, {\itshape Quaternionic 
   Plueker formula and Dirac eigenvalue estimates}, work in progress.
\bibitem{g}
   P. J. Giblin, {\itshape Graphs, Surfaces and Homology}, Chapman and 
   Hall Ltd, London, 1977.
\bibitem{k}
   N. Korevaar,
   {\itshape Upper bounds for eigenvalues of conformal metrics}, 
   J. Diff. Geom. 37, 73-93 (1993).
\bibitem{ly}
   P. Li and S-T. Yau,
   {\itshape Eigenvalues of a compact Riemannian manifold}, 
   Proc. of Symposia in Pure Math. 36, 205-239 (1980).
\bibitem{lnr}
   L. L. de Lima, V. F. de Sousa Neto and W. Rossman,
   {\itshape Lower bounds for index of Wente tori}, 
   to appear in Hiroshima Math. Journal.
\bibitem{lr}
     F. Lopez and A. Ros, 
     {\itshape Complete minimal surfaces with index one and stable 
     constant mean curvature surfaces}, 
     Comment. Math. Helvetici 64, 34-43 (1989).
\bibitem{j}
   C. Jaggy,
   {\itshape On the classification of constant mean curvature tori in 
   $\R^3$}, 
   Comment. Math. Helvetici 69, 640-658 (1994).
\bibitem{ps}
   U. Pinkall and I. Sterling,
   {\itshape On the classification of constant mean curvature tori}, 
   Annals of Math. 130, 407-451 (1989).
\bibitem{pw} 
   M. H. Protter and H. F. Weinberger, 
   {\itshape Maximum Principles in Differential Equations}, 
   Springer-Verlag (1984).  
\bibitem{rs}
   M. Reed and B. Simon,
   {\itshape IV Analysis of Operators}, 
   Academic Press, Inc. (1978).
\bibitem{r1}
   W. Rossman,
   {\itshape The Morse index of Wente tori}, 
   to appear in Geom. Dedicata.
\bibitem{r2}
   W. Rossman,
   {\itshape Wente tori and Morse index}, 
   An. Acad. Bras. Ci. 71, 607-613 (1999).
\bibitem{s}
     A. M. Silveira, 
     {\itshape Stability of complete noncompact surfaces 
     with constant mean curvature}, 
     Math. Ann. 277, 629-638 (1987).
\bibitem{u} 
     H. Urakawa, 
     {\itshape Geometry of Laplace-Beltrami operator 
     on a complete Riemannian manifold}, Advanced Studies in Pure 
     Mathematics, Progress in Differential Geometry 22, 347-406 (1993).  
\bibitem{wa}
     R. Walter, 
     {\itshape Explicit examples to the $H$-problem of Heinz Hopf}, 
     Geom. Dedicata 23, 187-213 (1987).
\bibitem{we}
     H. Wente, 
     {\itshape Counterexample of a conjecture of H. Hopf}, 
     Pacific J. Math. 121, 193-243 (1986).
\end{thebibliography}
\end{document}